\documentclass[a4paper,conference,10pt]{IEEEtran}

\usepackage{amssymb,amsmath,verbatim,bm,color,times,mathtools,multicol,relsize}
\usepackage{graphicx,tikz}
\usepackage{psfrag}
\usepackage[normalem]{ulem}
\usepackage{cite}
\usepackage{xcolor}
\definecolor{tom}{rgb}{0,0.5,0}
\usepackage{siunitx}
\usepackage{tabularx}
\usepackage{multirow,array,booktabs}
\usepackage{caption}
\usepackage{subcaption}
\usepackage{wrapfig}
\usepackage{textcomp}

\newcommand{\smallsum}{\mathop{\mathsmaller{\sum}}\nolimits}

\newcommand{\sumnl}{\sum\nolimits}

\renewcommand{\tabcolsep}{5pt}

\begin{document}
	
	\title{ \fontsize{16}{18}\selectfont Polynomial Chaos reformulation in Nonlinear Stochastic Optimal Control \\ with application on a drivetrain subject to bifurcation phenomena}
		
	\author{
		\IEEEauthorblockN{Tom Lefebvre}
		\IEEEauthorblockA{
			\textit{Department of Electrical Energy,} \\ \textit{Metals, Mechanical Constructions} \\ \textit{\& Systems; Ghent University}\\
			Ghent, Belgium \\
			\texttt{tom.lefebvre@ugent.be}
			}
		\and
		\IEEEauthorblockN{Frederik De Belie}
		\IEEEauthorblockA{
			\textit{Department of Electrical Energy,} \\ \textit{Metals, Mechanical Constructions} \\ \textit{\& Systems; Ghent University}\\
			Ghent, Belgium \\
			\texttt{frederik.debelie@ugent.be}
		}
		\and
		\IEEEauthorblockN{Guillaume Crevecoeur}
		\IEEEauthorblockA{
			\textit{Department of Electrical Energy,} \\ \textit{Metals, Mechanical Constructions} \\ \textit{\& Systems; Ghent University}\\
			Ghent, Belgium \\
			\texttt{guillaume.crevecoeur@ugent.be}
		}
		}
		
	\maketitle
	
	\begin{abstract}
		This paper discusses a method enabling optimal control of nonlinear systems that are subject to parametric uncertainty. A stochastic optimal tracking problem is formulated that can be expressed in function of the first two stochastic moments of the state. The proposed formulation allows to penalize system performance and system robustness independently. The use of polynomial chaos expansions is investigated to arrive at a computationally tractable formulation expressing the stochastic moments in function of the polynomial expansion coefficients rigorously. It is then demonstrated how the stochastic optimal control problem can be reformulated as a deterministic optimal control problem in function of these coefficients. The proposed method is applied to find a robust control input for the start-up of an eccentrically loaded drive train that is inherently prone to bifurcation behaviour. A reference trajectory is chosen to deliberately provoke a bifurcation. The proposed framework is able to avoid the bifurcation behaviour regardlessly. 

	\end{abstract}
	
	
	\section{Introduction}
	The Optimal Control (OC) method is consolidating into an affordable mathematical instrument offering a decisive solution to various challenges in a wide range of process control applications.
	Successful application of the method is well documented in the control literature, both in terms of preoperational planning, \cite{richter2016polynomial,diehl2006fast}, as well as to realise real-time state feedback control policies, \cite{torchio2015real,klanvcar2007tracking}. Under deterministic conditions, where only the exogenous disturbances are stochastic of nature, excellent performance is obtained. 
	
	Its success is greatly owing to the predictive capabilities of the nonlinear system model at hand. Reality urges one to reconsider this deterministic modelling assumption however, since model uncertainty may be inherent to many applications, jeopardizing the very predictive capability of the model; e.g. parametric uncertainty, model-plant mismatch, etc. 
	
	Techniques that account for parametric uncertainty in OC have become collected under the heading stochastic OC. In stochastic OC, one aims to find a solution that minimizes the cost in expectation or with some given probability \cite{fisher2011optimal,mesbah2014stochastic}; that as a less conservative alternative for robust OC where one is interested in worst-case solutions only \cite{cannon2011stochastic,beyer2007robust}. 
	
	Contrary to reformulations that address the expected value of the cost functional \cite{mesbah2014stochastic,YuningJiang2017}, we propose a specific reformulation that takes into account the robustness of the computed control signal as well. To that end we add the Frobenius norm of the state covariance to the cost kernel in concordance to its use in static robust optimization \cite{MANDUR2012549}. The performance and robustness are penalized independently as a result. It will show that this formulation allows to express the OC problem in function of the stochastic moments of the state.
	
	The key challenge in stochastic OC can thus be identified as the propagation of the parametric uncertainty through the system model. In the case of nonlinear systems, this is nontrivial and often an analytical expression is unlikely, if not impossible to be found, so that we are forced to rely on numerical approximation strategies to quantify the moments. 
	
	An ubiquitous method to do so is Monte Carlo (MC) sampling, or derivatives thereof, such as the Sigma Point approach \cite{YuningJiang2017,blackmore2010probabilistic}. 
	Such approaches are curtailed by the computational challenge that comes with the numerous forward simulations required to achieve an acceptable degree of accuracy or are limited to symmetrical distributions.
	
	The generalised Polynomial Chaos (gPC) framework provides an advantageous setting to quantify and propagate parametric uncertainty and is computationally more efficient over MC approaches \cite{kim2013wiener}. The propagation of input uncertainty to the output variables is realised by modelling the explicit nonlinear output function through a polynomial expansion. By choosing the polynomial basis so that it is orthogonal with respect to the probability density function of the uncertain parameters, a strong theoretical relation is established that allows to express the statistical moments as a function of the polynomial coefficients. The stochastic OC can then be reformulated as a deterministic OC problem.
	
	In the current paper we propose a stochastic OC problem class that can be expressed in function of the first two stochastic moments of the state. We discuss the gPC framework as a tool to numerically approximate the statistical moments in function of the polynomial expansion coefficients rigorously. 
	
	The computational framework is applied on the start-up behaviour of an eccentrically loaded and thus nonlinear drivetrain that is prone to bifurcation behaviour, i.e. a sudden qualitative change in dynamic behaviour when a system parameter is changed. We
	consider parametric uncertainty that affects both the initial state as the system dynamics, thus creating an opportunity to handle both. 

	\section{Problem formulation}
	\label{sec:stochastic}
	We consider stochastic OC problems of the form
	\begin{equation}
	\label{eq:problem}
	\begin{aligned}
	&\min_{\mathbf{u}(\cdot)} ~ \Omega\left[ J(\mathbf{x}(\cdot,\boldsymbol{\omega}),\mathbf{u}(\cdot))\right] \\
	\text{s.t. } &\left\lbrace\begin{aligned}
	\mathbf{x}(0,\boldsymbol{\omega}) &= \mathbf{h}(\boldsymbol{\omega}), \\
	\dot{\mathbf{x}}(t,\boldsymbol{\omega}) &= \mathbf{f}(t,\mathbf{x}(t,\boldsymbol{\omega}),\mathbf{u}(t),\boldsymbol{\omega}), &t\in\mathcal{T}\\
	\uline{\mathbf{u}}&\leq\mathbf{u}(t)\leq\overline{\mathbf{u}}, &t\in\mathcal{T}
	\end{aligned}\right.
	\end{aligned}
	\end{equation}
	
	The random vector $\boldsymbol{\omega}\in\Gamma\subseteq\mathbb{R}^{n_\omega}$ represents an uncertain time-invariant parameter. The parametric state variable is denoted by $\mathbf{x}:\mathbb{R}\times\mathbb{R}^{n_\omega}\rightarrow\mathbb{R}^{n_x}$. The vector $\mathbf{u}:\mathbb{R}\rightarrow\mathbb{R}^{n_u}$ represents a deterministic control signal bounded by $\underline{\mathbf{u}}$ and $\overline{\mathbf{u}}$. We assume that the entries of $\boldsymbol{\omega}$ are independent random variables, $\omega_k$, with known probability distribution functions, $\rho_k:\mathbb{R}\rightarrow\mathbb{R}_{\geq0}$. 
	Functions $\mathbf{f}:\mathbb{R}\times\mathbb{R}^{n_x}\times\mathbb{R}^{n_u}\times\mathbb{R}^{n_\omega}\rightarrow\mathbb{R}^{n_x}$ and $\mathbf{h}:\mathbb{R}^{n_\omega}\rightarrow\mathbb{R}^{n_x}$ represent the governing time-variant dynamic system equation and the initial state of the system, respectively. The dependency of both $\mathbf{f}$ and $\mathbf{h}$ on the random vector, $\boldsymbol{\omega}$, exemplifies that the uncertainty can root both from the initial state as from the system model. 
	
	A stochastic formulation of the control objective is obtained through the probabilistic operator, $\Omega$, that associates a single objective value to the entire stochastic space $\Gamma$. 
	
	For the underlying deterministic optimal control objective, we take interest in the quadratic cost functional
	\begin{equation}
	\label{eq:costfunctional}
	J \doteq \int_{\mathcal{T}} \left\|\mathbf{x}(\tau,\boldsymbol{\omega})-\mathbf{r}(\tau)\right\|_{{\mathrm{Q}}}^2+ \|\mathbf{u}(\tau)\|_{{\mathrm{R}}}^2\text{d}\tau
	\end{equation}
	where $\left\|\mathbf{z}\right\|^2_{{\mathrm{A}}}\doteq\mathbf{z}^\top{\mathrm{A}}\mathbf{z}$ and for a given time-dependent reference trajectory, $\mathbf{r}:\mathbb{R}\rightarrow\mathbb{R}^{n_x}$, and, positive definite weight matrices ${\mathrm{Q}}\in\mathbb{R}^{n_x\times n_x}$ and ${\mathrm{R}}\in\mathbb{R}^{n_u\times n_u}$.
	
	In order to provide a tractable stochastic reformulation of the deterministic control objective, we define the probabilistic operator, $\Omega$, as the weighted sum of the expected cost value and the integrated Frobenius norm of the state covariance matrix. This definition provides a trade-off between maximized performance, by the expected cost value, and minimized uncertainty, i.e. robustness, of the final trajectory. We deliberately avoid to penalize the variability of the tracking error instead of that of the optimized trajectory itself, so that performance and robustness of the solution can be penalized independently by altering the factor $0<\epsilon\leq1$. 
	\begin{equation}
	\label{eq:stochasticoperator}
	\begin{aligned}
	\Omega[J] &\doteq \epsilon\mathbb{E}[J] + (1-\epsilon)\int_\mathcal{T} \|\mathrm{cov}[\mathbf{x}]\|_F^2\text{d}\tau \\
	&\begin{multlined}=
	\epsilon \int_\mathcal{T} \mathbb{E}[\|\mathbf{x}\|_\mathrm{Q}^2]+\left\|\mathbf{r}\right\|_{{\mathrm{Q}}}^2-2\mathbf{r}^\top{\mathrm{Q}}\mathbb{E}\left[\mathbf{x}\right]+\left\|\mathbf{u}\right\|^2_{{\mathrm{R}}}\text{d}\tau \\
	+ (1-\epsilon) \int_\mathcal{T} \|\mathbb{E}[\mathbf{x}]\|^2 - \mathbb{E}[\|\mathbf{x}\|^2]\text{d}\tau
	\end{multlined}
	\end{aligned}
	\end{equation}
	
	Under these assumptions numerical solution of (\ref{eq:problem}) will require a means to evaluate the first two stochastic moments of the state. Classically these statistical moments would be obtained by MC methods. The number of forward system simulation would grow prohibitively large however considering the iterative solution of (\ref{eq:problem}). The gPC framework is known to outperform MC methods and is advanced alternatively. 
	

	\section{Polynomial Chaos Expansions}
	\subsection{Wiener-Askey polynomial chaos}
	According the gPC framework \cite{xiu2002wiener}, a sufficiently smooth function $x:\mathbb{R}^{n_\omega}\rightarrow\mathbb{R}$, can be modelled as an infinite summation of polynomials. From a computational perspective we take interest in the $d$-th order approximation. That is, let $\mathcal{P}^d_{n_\omega}$ be the $n_\omega$-variate polynomial space of at most degree $d$ and $\boldsymbol{\Psi} = \{\Psi_1,\cdots,\Psi_p\}$\footnote{$\boldsymbol{\Psi}$ can be generated from the univariate bases $\boldsymbol{\Phi}^{(k)} = \{\Phi^{(k)}_0,\dots,\Phi^{(k)}_{d}\}$. Consider $\Psi_{|\mathbf{i}|\leq d}(\boldsymbol{\omega}) = \prod_{k=1}^{n_\omega}\Phi^{(k)}_{i_k}(\omega_k)$, where $\mathbf{i} = (i_1,\dots,i_l)$ and $|\mathbf{i}|=\sum_{k=1}^{n_\omega} i_k$. For notational convenience, we exploit the existing bijection between multi-index $\mathbf{i}$ and the index $i$ taking values in $\{1,\dots,\binom{n_\omega+d}{n_\omega}\}$.} 
	be a basis of $\mathcal{P}^d_{n_\omega}$ with $p=\binom{n_\omega+d}{n_\omega}$.
	
	The $d$-th order polynomial approximation $x_{(d)}(\boldsymbol{\omega})$ of $x(\boldsymbol{\omega})$ is then given by (\ref{eq:pcapprox}) for given polynomial coefficients $\tilde{x}_i$.
	\begin{equation}
	\label{eq:pcapprox}
	x_{(d)}(\boldsymbol{\omega}) = \sum_{i=1}^{p} \tilde{x}_i \Psi_i(\boldsymbol{\omega}) \xrightarrow{d\rightarrow\infty} x(\boldsymbol{\omega}) = \sum_{i=1}^{\infty} \tilde{x}_i \Psi_i(\boldsymbol{\omega}) 
	\end{equation}
	
	Within the context of uncertainty propagation this model allows to establish advantageous computational conditions by a conscious choice of the basis, $\boldsymbol{\Psi}$. First, let $\boldsymbol{\omega}$ be composed of $n_\omega$ independently distributed random variables, $\omega_k$, with known supports and power density functions (PDF) $\rho_k:\Gamma_k\subseteq\mathbb{R}\rightarrow\mathbb{R}_{\geq0}$. The joint support, $\Gamma$, and PDF, $\rho$, are hence given by $\Gamma=\bigotimes_{k}\Gamma_k$ and $\rho=\prod_{k}\rho_k$. We are interested in propagating the input uncertainty on $\boldsymbol{\omega}$ to the output $x$. By choosing the univariate bases, $\boldsymbol{\Phi}^{(k)}$, so that they satisfy an orthogonality condition with respect to the PDF associated to their respective variable, the statistical moments can be obtained in function of the polynomial coefficients. 
	
	In the context of polynomial approximations, orthogonality of a basis, $\boldsymbol{\Psi}$, is defined indirectly through the inner product, $\langle\Psi_i,\Psi_j\rangle \doteq \int_{\Gamma} \Psi_i\Psi_j\rho\text{d}\boldsymbol{\omega}$. Basis $\boldsymbol{\Psi}$ is said to be orthogonal for given $\rho$ if $\langle\Psi_i,\Psi_j\rangle = \delta_{ij}\langle\Psi_i^2\rangle$. The inner product determines a projection operator on the polynomial space, $\mathcal{P}^d_{n_\omega}$, and is arbitrarily defined by the weighting, $\rho$. 
	
	Interestingly the $n$-th stochastic moment, $\mu'_n$, of $x(\boldsymbol{\omega})$ can then be approximated in function of the coefficients, $\tilde{x}_i$
	\begin{equation}
	\label{eq:moments}
	\mu'_n \approx \mathbb{E}\left[x_{(d)}^n\right] = \sumnl_{i_1} \cdots \sumnl_{i_n} \tilde{x}_{i_1} \cdots \tilde{x}_{i_n} \langle \Psi_{i_1} \cdots \Psi_{i_n} \rangle
	\end{equation}
	
	
	Now let’s choose each univariate basis $\boldsymbol{\Phi}^{(k)}$ so that it satisfies the orthogonality condition with respect to the PDF, $\rho_k$, implying that $\langle\Psi_i,\Psi_j\rangle =  \prod_k\langle\Phi^{(k)}_{i_k},\Phi^{(k)}_{j_k}\rangle$ due to the independency of the $\omega_k$s. Common probability measures are associated to a polynomial basis by the Wiener-Askey scheme \cite{xiu2007efficient}, for arbitrary PDF a basis can be generated. 
	
	Starting from (\ref{eq:moments}) it is then easily verified that the first two stochastic moments, $\mu'_1$, and, $\mu'_2$, reduce to
	\begin{align}
	\label{eq:2moments}
	\mu'_1 = \tilde{x}_1\langle\Psi_1^2\rangle & & \mu'_2 = \sumnl_i \tilde{x}^2_i \langle\Psi_i^2\rangle
	\end{align}
	
	
	In conclusion, we can rewrite the approximation in a dense matrix form and generalise it to $n$-dimensional output models. The basis and coefficient vectors are defined as $\boldsymbol{\Psi} = [\Psi_1,\dots,\Psi_p]\in\mathbb{R}^p$ and ${\tilde{\mathbf{X}}} = [\mathbf{\tilde{x}}_1;\dots;\mathbf{\tilde{x}}_p]\in\mathbb{R}^{n p}$.
	\begin{align}
	\mathbf{x}_{(d)}(\boldsymbol{\omega}) = \boldsymbol{\Psi}(\boldsymbol{\omega}) \otimes {\mathrm{I}} \cdot {\tilde{\mathbf{X}}}, & & \mathbf{x}:\mathbb{R}^{n_\omega}\rightarrow\mathbb{R}^{n}
	\end{align}
	
	\subsection{Coefficient determination}
	The framework as such accounts for a computationally advantageous mathematical setting allowing to quantify the output uncertainty in function of the polynomial coefficients rigorously. In this section we provide an overview of three common methods in the literature to determine the polynomial coefficients \cite{xiu2007efficient,hadigol2017least,cheng2010collocation}. We also point out that under certain conditions two of these methods are numerically equivalent. 
	
	All of the discussed techniques are non-intrusive methods, i.e. they retrieve the desired coefficients by means of a set of deterministic evaluations of the model response corresponding an input set, $\left\{\boldsymbol{\omega}_j\right\}_{j=1}^q$. The term non-intrusive indicates that these methods can be applied using the deterministic code associated with the forward model $\mathbf{x}$ without modification. From that perspective these techniques are similar to the MC approach. The fundamental difference is in that the estimation of the statistical moments is realized through a mathematical detour. First the forward model is approximated by the polynomial series. By proper definition of the polynomial basis the orthogonality property allows to calculate the statistical moments exactly and directly from the polynomial coefficients. The focus of approximation in the gPC framework is therefore on modelling the response function whilst that of the MC approach is on the direct estimation of the statistical moments. Therefore the level of accuracy depends on the capability of the basis to capture the nonlinearity of the forward model rather than on the capacity of the sample generation algorithm to properly represent the input uncertainty by the spatial distribution of the sample points. It has been documented frequently that an equivalent level of accuracy can be achieved with only a fraction of the input points of any MC approach \cite{7130673,6600994,kewlani2012polynomial}.
	
	The methods are explicated with respect to the function $\mathbf{x}:\mathbb{R}^{n_\omega}\rightarrow\mathbb{R}^{n}$.
	\paragraph{Projection method} The projection method (PM) projects the response function $\mathbf{x}$ onto the polynomial space $\mathcal{P}^d_{n_\omega}$ by application of the operator $\langle\cdot,\Psi_i\rangle$. The projection coefficients are then given by (\ref{eq:galproj-exact}). For general nonlinear models, this integral cannot be determined explicitly and is approximated by a Gauss-quadrature rule (\ref{eq:galproj-quad}). A quadrature is defined by a set of $q$ collocation points and corresponding weights, $\{(\boldsymbol{\omega}_j,w_j)\}_{j=1}^q$. For additional details on quadrature rules and collocation sets we refer to, e.g. \cite{hadigol2017least}.
	\begin{subequations}
		\begin{align}
		\label{eq:galproj-exact}
		\mathbf{\tilde{x}}_i \langle\Psi_i^2\rangle = \left\langle \mathbf{x},\Psi_i\right\rangle &= \int_{\Gamma} \mathbf{x}\cdot\Psi_i\rho\text{d}\boldsymbol{\omega} \\
		\label{eq:galproj-quad}
		\left\langle \mathbf{x},\Psi_i\right\rangle &\approx \sumnl_{j} \mathbf{x}(\boldsymbol{\omega}_j)\Psi_i(\boldsymbol{\omega}_j) w_j
		\end{align}
	\end{subequations}
	
	This operation can be represented compactly by an affine transformation between the collocation and the coefficient vector spaces
	\begin{equation}
	{\tilde{\mathbf{X}}} = {\mathrm{D}}^{-1}{\mathrm{\Psi}}^\top {{\mathrm{W}}}\otimes {\mathrm{I}}\cdot{\hat{\mathbf{X}}}
	\end{equation}
	where the matrices are defined as $\mathrm{D}_{ij} = \langle\Psi_i,\Psi_j\rangle$, ${\mathrm{W}}_{ii} = w_i$ and $\Psi_{ij} = \Psi_j(\boldsymbol{\omega}_i)$ respectively and introducing the collocation vector, $\hat{{\mathbf{X}}} = [\mathbf{x}(\boldsymbol{\omega}_1);\dots;\mathbf{x}(\boldsymbol{\omega}_q)]\in\mathbb{R}^{nq}$.
	
	\paragraph{Least-Squares} Another common approach is to perform a classical regression analysis \cite{hadigol2017least}.  Given that again $q$ collocation points $\{\mathbf{x}(\boldsymbol{\omega}_j)\}_{j=1}^q$ are available, the coefficients are determined by the solution of the least-squares (LS) problem (here defined for a 1D output)
	\begin{equation}
	\label{eq:LScoeff}
	\min_{{\mathrm{\tilde{X}}}} \|{\hat{\mathbf{X}}} - {\Psi}\cdot{\tilde{\mathbf{X}}}\|^2_2
	\end{equation}
	
	We remark that in this case the choice of the collocation points is more arbitrary. Typically a collocation set is obtained through design of computer experiments. Examples are Hammersley sampling, Latin Hypercube sampling, D-optimal sampling etc. \cite{hadigol2017least}. The solution of (\ref{eq:LScoeff}) can again be represented compactly by an affine transformation between the collocation and the coefficient vector spaces
	\begin{equation}
	{\tilde{\mathbf{X}}} =({\Psi}^\top{\Psi})^{-1}{\Psi}^\top\otimes {\mathrm{I}} \cdot{\hat{\mathbf{X}}}
	\end{equation}
	
	We remark that the LS approach differs from the projection method in this aspect that the projection only recollects the first $p$ terms of the infinite summation in (\ref{eq:pcapprox}), whilst the LS approach tries to compensate for the truncated tail.
	
	\paragraph{Generalised Least-Squares} The last method entails a generalisation of the LS approach \cite{cheng2010collocation}. Conceptually, the coefficients are determined so that they minimize the expected squared error between the true function and its polynomial approximation. When the occurring integral form is approximated numerically, explicit expression of the expected value coincides with that of a generalised LS (gLS) formulation
	\begin{equation}
	\begin{aligned}
	\label{eq:gLS}
	\min_{{\tilde{\mathrm{X}}}} \mathbb{E}[\|\mathbf{x} - \mathbf{x}_{(d)}\|_2^2] &= \int_{\Gamma} \|\mathbf{x} - \mathbf{x}_{(d)}\|_2^2\rho\text{d}\underline{\omega} \\
	& \approx \sumnl_{j} \|\mathbf{x}(\boldsymbol{\omega}_j) - {\mathbf{x}}_{(d)}(\boldsymbol{\omega}_j)\|_2^2w_j \\
	& = \|{\hat{\mathbf{X}}} - {\Psi}\cdot{\tilde{\mathbf{X}}}\|^2_{\mathrm{W}}
	\end{aligned}
	\end{equation}
	
	The corresponding matrix transformation is then given by
	\begin{equation}
	{\tilde{\mathbf{X}}} =({\Psi}^\top{\mathrm{W}}{\Psi})^{-1}{\Psi}^\top{\mathrm{W}}\otimes {\mathrm{I}} \cdot{\hat{\mathbf{X}}}
	\end{equation}
	
	This and the projection method are equivalent in the case that a Gaussian-quadrature rule is adopted that is exact for polynomials of degree $2d$. Consider therefore that, if this is the case
	\begin{equation}
	\begin{aligned}
	({\Psi}^\top{\mathrm{W}}{\Psi})_{ij} &= \sumnl_k \Psi_i(\boldsymbol{\omega}_k)\Psi_j(\boldsymbol{\omega}_k)w_k \\
	&= \langle\Psi_i,\Psi_j\rangle = \mathrm{D}_{ij}
	\end{aligned}
	\end{equation}
	
	This shows that the projection and the gLS method are both formally and numerically equivalent. The formal equivalence of (\ref{eq:galproj-exact}) and (\ref{eq:gLS}) is a fundamental result in polynomial function approximation, see for example \cite{cameron1947orthogonal}.
	
	\subsection{Propagation of uncertainty in dynamical systems}
	As apparent in the previous section, any nonintrusive method can be represented by an affine mapping from the collocation to the coefficient vector space. For time-variant output models, consider the following generalisation
	\begin{equation}
	\label{eq:coefficienttransform}
	{\tilde{\mathbf{X}}}(t) = {\mathrm{A}}\cdot {\hat{\mathbf{X}}}(t)
	\end{equation}
	
	Starting hereof, we introduce two particularly different coefficient determination methods for when the model,  $\mathbf{x}(t,\boldsymbol{\omega})$, satisfies the nonlinear ordinary differential equation (ODE)
	\begin{equation}
	\left\lbrace
	\begin{aligned}
	\mathbf{x}(0,\boldsymbol{\omega}) &= \mathbf{h}(\boldsymbol{\omega}) \\
	\dot{\mathbf{x}}(t,\boldsymbol{\omega}) &= \mathbf{f}(t,\mathbf{x}(t,\boldsymbol{\omega}),\boldsymbol{\omega}), & t\in\mathcal{T}
	\end{aligned}
	\right.
	\end{equation}
	
	\paragraph{Decoupled coefficient determination} The decoupled approach is a direct application of (\ref{eq:coefficienttransform}). The entries of the time-variant collocation vector are simply determined by considering $q$ independent ODEs, formally represented as
	\begin{equation}
	\label{eq:collocationBVP}
	\left\lbrace\begin{aligned}
	{\hat{\mathbf{X}}}(0) &= \begin{bmatrix}
	\mathbf{h}(\boldsymbol{\omega}_1) \\
	\vdots \\
	\mathbf{h}(\boldsymbol{\omega}_q) 
	\end{bmatrix} \doteq \hat{{\mathbf{H}}} \\
	\dot{{\hat{\mathbf{X}}}}(t) &= \begin{bmatrix}
	\mathbf{f}(t,\mathbf{x}(t,\boldsymbol{\omega}_1)),\boldsymbol{\omega}_1) \\ \vdots \\ \mathbf{f}(t,\mathbf{x}(t,\boldsymbol{\omega}_q)),\boldsymbol{\omega}_q)
	\end{bmatrix} \doteq \hat{{\mathbf{F}}}(t,{\hat{\mathbf{X}}}(t))
	\end{aligned}\right.
	\end{equation}
	
	${\tilde{\mathbf{X}}}(t)$ can then be determined by application of (\ref{eq:coefficienttransform}).
	
	\paragraph{Coupled coefficient determination} The coupled approach considers the time derivative of (\ref{eq:coefficienttransform}) and replaces the state signal by its polynomial approximation. That is
	\begin{equation}
	\begin{aligned}
	\dot{{\tilde{\mathbf{X}}}}(t) &= {\mathrm{A}}\cdot \dot{{\hat{\mathbf{X}}}}(t) = {\mathrm{A}}\cdot\begin{bmatrix}
	\mathbf{f}(t,\mathbf{x}(t,\boldsymbol{\omega}_1),\boldsymbol{\omega}_1) \\ \vdots \\ \mathbf{f}(t,\mathbf{x}(t,\boldsymbol{\omega}_q),\boldsymbol{\omega}_q) 
	\end{bmatrix} \\
	&\approx {\mathrm{A}}\cdot\begin{bmatrix}
	\mathbf{f}(t,\tilde{\boldsymbol{\Psi}}(\boldsymbol{\omega}_1) \cdot {\tilde{\mathbf{X}}}(t),\boldsymbol{\omega}_1) \\ \vdots \\ \mathbf{f}(t,\tilde{\boldsymbol{\Psi}}(\boldsymbol{\omega}_q) \cdot {\tilde{\mathbf{X}}}(t),\boldsymbol{\omega}_q)
	\end{bmatrix} \doteq \tilde{{\mathbf{F}}}(t,{\tilde{\mathbf{X}}}(t))
	\end{aligned}
	\end{equation}
	where $\boldsymbol{\tilde{\Psi}} = \boldsymbol{\Psi}\otimes{\mathrm{I}}$. 
	
	Hence we obtain an ODE in function of the coefficients\footnote{Note that the decoupled and coupled method are equivalent when $p=q$ due to the uniqueness of the inverse of $\mathrm{A}$.}.
	\begin{equation}
	\label{eq:coefficientBVP}
	\left\lbrace\begin{aligned}
	\tilde{{\mathbf{X}}}(0) &= {\mathrm{A}}\cdot {\hat{\mathbf{H}}} \\
	\dot{{\tilde{\mathbf{X}}}}(t) & =  {\mathbf{F}}(t,{\tilde{\mathbf{X}}}(t))
	\end{aligned}\right.
	\end{equation}
	


	\section{Stochastic Optimal Control with gPC}
	\label{sec:stochastic_gPC}
	Here we illustrate how the gPC expansion can be of use in the stochastic OC framework discussed in section \ref{sec:stochastic}. 
	
	We introduce the piecewise linear control policy, $\mathbf{u}:\mathbb{R}\times\mathbb{R}^{n_u (n_t+1)}\rightarrow\mathbb{R}^{n_u}$ parameterised by $\mathbf{\mathbf{U}} = [\mathbf{u}_0;\dots;\mathbf{u}_{n_t}]\in\mathbb{R}^{n_u (n_t+1)}$. Here, $d^k$, is the discrete delay operator,  $\mathrm{tri}(t)=\max(0,1-|t|)$, the triangular function and $\Delta = \tfrac{\mathcal{T}}{n_t}$.
	\begin{align}
	\mathbf{u}(t;{\mathbf{U}}) = \sumnl_{k=0}^{n_t} \mathbf{u}_k d^{k}\mathrm{tri}\left(\tfrac{t}{\Delta}\right), & & t\in\mathcal{T}
	\end{align}
	
	We further assume that the time dependency of the system dynamics results solely from the time-variant input policy. The dynamics of a controlled system can as such be represented formally as given by
	\begin{equation}
	\dot{\mathbf{x}}(t,\boldsymbol{\omega}) = \mathbf{f}(t,\mathbf{x}(t,\boldsymbol{\omega})) \doteq \mathbf{f}(\mathbf{x}(t,\boldsymbol{\omega}),\mathbf{u}(t)) 
	\end{equation}
	
	\subsection{Polynomial Chaos reformulation}
	We define the stochastic cost functional $K$ as $ K(\mathbf{x},\mathbf{u}) = \Omega[J(\mathbf{x},\mathbf{u})]$. As a result of the specific definition of $J$ (\ref{eq:costfunctional}) and $\Omega$ (\ref{eq:stochasticoperator}), functional $K$ can be expressed as a quadratic cost in function of the coefficient, ${\tilde{\mathbf{X}}}$, and control vector, ${\mathbf{U}}$.
	
	The expected cost value reduces to (see Appendix \ref{appendix:details_PCeSOC})
	\begin{equation}
	\mathbb{E}[J] = \int_\mathcal{T} \|{\tilde{\mathbf{X}}}-{{\tilde{\mathbf{R}}}}\|^2_{\mathrm{D}\otimes\mathrm{Q}}\text{d}\tau + \|{\mathbf{U}}\|^2_{\mathrm{M}\otimes\mathrm{R}}
	\end{equation}
	where $\tilde{{\mathbf{R}}}$ is defined as $\mathbf{e}_1\otimes\mathbf{r}$, $\mathbf{e}_1\in\mathbb{R}^p$ and where ${\mathrm{M}}\in\mathbb{R}^{(n_t+1)\times(n_t+1)}$ is a tridiagonal matrix, related to the evaluation of the control contribution in (\ref{eq:costfunctional}), that for this particular parameterization of the control policy, $\mathbf{u}(t)$, is equal to
	\begin{equation}
	{\mathrm{M}}= \frac{\Delta}{6}\left[\begin{smallmatrix}
	2 & 1 & & & \\
	1 & 4 & 1& & \\
	&\ddots & \ddots & \ddots& \\
	& & 1 & 4 & 1 \\
	& & & 1 & 2 
	\end{smallmatrix}\right] 
	\end{equation} 
	
	The Frobenius norm of the state covariance reduces to
	\begin{equation}
	\int_{\mathcal{T}}\|\mathrm{cov}[{\mathbf{x}}]\|_F^2 \text{d}t \approx  \int_{\mathcal{T}}\|{\tilde{\mathbf{X}}}\|^2_{\mathrm{E}\otimes\mathrm{I}} \text{d}t
	\end{equation}
	where ${\mathrm{E}} = {\mathrm{D}}-\mathbf{e}_1\mathbf{e}_1^\top$ (see Appendix \ref{appendix:details_PCeSOC}).
	
	The results obtained heretofore can be summarized elegantly into the following compact expression for $K$
	\begin{multline}
	\label{eq:costfunctional_stochastic} K(\mathbf{x},\mathbf{u}) \approx 
	K({\tilde{\mathbf{X}}},{\mathbf{U}}) \doteq K(\mathbf{x}_{(d)},\mathbf{u})  = \\ \int_\mathcal{T} \epsilon\|{\tilde{\mathbf{X}}}-{\tilde{\mathbf{R}}}\|^2_{\mathrm{D}\otimes\mathrm{Q}}+(1-\epsilon) \|{\tilde{\mathbf{X}}}\|^2_{\mathrm{E}\otimes\mathrm{I}}\text{d}\tau + \epsilon\|{\mathbf{U}}\|^2_{\mathrm{M}\otimes\mathrm{R}} 
	\end{multline}
	
	Considering, (\ref{eq:coefficienttransform}) and (\ref{eq:collocationBVP}), either (\ref{eq:coefficientBVP}), and the approximation of $K$ in (\ref{eq:costfunctional_stochastic}), problem (\ref{eq:problem}) can be reformulated as 
	\begin{equation}
	\label{eq:stochasticproblem}
	\begin{aligned}
	&\min_{\uline{\mathbf{U}}\leq{\mathbf{U}}\leq\overline{\mathbf{U}}} ~ K({\tilde{\mathbf{X}}},{\mathbf{U}}) \\
	\text{s.t. } & (\ref{eq:coefficienttransform})\text{ and }(\ref{eq:collocationBVP}) \text{, or }(\ref{eq:coefficientBVP})
	\end{aligned}
	\end{equation}
	
	The former is a deterministic optimal control problem in function of the coefficients and can be solved accordingly.
	
	\subsection{Remarks}
	Regarding the formulation in (\ref{eq:problem}) it is possible to consider $\mathbf{r}$ as a control variable. Such would allow to find a reference trajectory, $\mathbf{r}^*$, that is by construction easy to track under given uncertainty of the system parameters. Regarding the formulation in (\ref{eq:costfunctional_stochastic}), $\mathbf{r}^*$ should then simply coincide with $\tilde{\mathbf{x}}^*_1$. 
	
	The problem is thus equivalent to the current formulation if extended with a(n) (equality) constraint on the (expected) end state, e.g. $\tilde{\mathbf{x}}_1(T) = \mathbf{x}_T$, and with $\mathbf{r}=\tilde{\mathbf{x}}_1$. Notice that the integrand in $K$ then reduces to $\|{\tilde{\mathbf{X}}}\|^2_{\mathrm{E}\otimes\epsilon\mathrm{Q} +(1-\epsilon)\mathrm{I}}$.
	
	
	\section{Stochastic Optimal Control \\ of a Small Drivetrain}
	The framework elaborated above is applied to the start-up of a small drivetrain. The drivetrain consists of a flywheel with inertia $J=1$ and an eccentrically loaded shaft end. The shaft is radially loaded at a distance $r=1$ from its centre by a spring ($k = 1$) and damper element ($b = 0.5$) that are both suspended at a distance $l=1.5$ from the shaft centre. The nonlinear torque contributions resulting the eccentrically placed spring damper element are given by $T_k$ and $T_b$. 
	\newpage
	
	\begin{figure}[t!]
		\centering
		\includegraphics[trim=0cm 0cm 0cm 0cm,clip=true,width=.6\columnwidth]{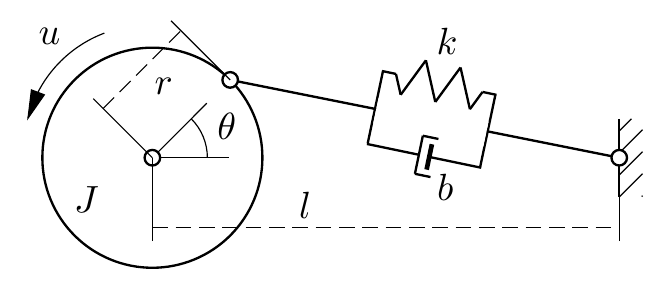}
		\caption{Mechanical system diagram.}
			\vspace*{-15pt}
	\end{figure}

	Note that $T_k$ is a function of the variable parameter $\theta_0$ that determines the axis rotation at rest as determined by the actual rest length of the spring element, $l_0$.
	\begin{subequations}
		\begin{align}
		T_k(\theta) &= k\left(1-\frac{\sqrt{r^2+l^2-2rl\cos\theta}}{\sqrt{r^2+l^2-2rl\cos\theta_0}}\right)rl\sin\theta \\
		T_b(\theta,\dot{\theta}) &= b \frac{r^2l^2\sin\theta^2}{r^2+l^2-2rl\cos\theta}\dot{\theta}
		\end{align}
	\end{subequations}
	
	The systems dynamics are then determined by the following state-space equation
	\begin{equation}
	J\ddot{\theta} + T_k(\theta) + T_b(\theta,\dot{\theta})= u 
	\end{equation}
	
	We consider stochasticity of $\theta_0=\tfrac{\pi}{2}+\omega\tfrac{\pi}{4},\omega\sim\mathcal{U}(-1,1)$, affecting start-up behaviour from rest, such that the uncertainty will contaminate both the dynamics as the initial state.
	
	\subsection{gPC versus MC for a step response}
	We first examine the modelling and uncertainty propagation capabilities of the gPC expansion framework, before venturing into the robustified control reformulation. The step response behaviour of the system is considered. Two characteristic responses can be distinguished. Either the input torque will be sufficient to obtain rotational motion, otherwise the axle will settle for an equilibrium without rotation. The uncertainty of the parameter $\theta_0$ thus implies that, for a certain range of step inputs, both regimes can be expected. Hence, the system will be prone to bifurcation behaviour.
	
	Two step inputs are considered, $u_1(t) = \tfrac{1}{2}\mathbb{H}(t)$ and $u_2(t) = 1\mathbb{H}(t)$. For $u_1$, the bifurcation behaviour will be present. With $u_2$, rotational motion is obtained regardless the value of $\theta_0$. We refer to these settings as \textit{scenario} $1$ and \textit{scenario} $2$, respectively. For every scenario, the system is simulated over a period of $10$ \si{\second}. The states are defined as $x_1 = \theta$ and $x_2 = \dot{\theta}$. The uniform distribution corresponds the Legendre polynomials \cite{xiu2007efficient}. We used Gaussian quadrature to obtain the collocation points. 
	
	Table \ref{tab:rmse} compares the root mean squared error (RMSE) between the surfaces, $x_{(d),1}(t,\omega)$ and $x(t,\omega)$. Surface $x_{(d),1}(t,\omega)$ is obtained with gPC for varying polynomial degree, $d$, and number of collocation nodes, $q$. Surface $x(t,\omega)$ is obtained with $500$ Monte Carlo (MC) simulations. The gPC surface is obtained with both the LS and the PM, as well as the coupled and decoupled coefficient determination methods. We present only the RMSE value obtained for the first state, results for the second state are similar.
	
	One can observe quasi exponential convergence concerning the polynomial degree, $d$. The value of $q$ does not have a very pronounced effect (as long as $q\geq p = d+1$), especially when the polynomial degree increases. The difference between the PM and the LS approach is close to negligible, although the PM approach performs best in any case. Also the performance of the decoupled and coupled dynamic coefficient determination are quite comparable. Nonetheless a slight preference is observed for the coupled approach w.r.t. scenario $1$ and for the decoupled approach w.r.t. scenario $2$. This might be explained considering that the decoupled approach is more off a post-processing technique and that if the distribution of collocation points is too coarse, the discontinuity is under sampled and cannot be modelled exactly, which is less the case for the coupled approach that solves the problem in the coefficient domain.
	
\begin{figure}[t!]
		\centering
		\begin{subfigure}[b]{0.43\columnwidth}
			\caption{$x_1(t,\omega)$}
			\vspace{-6pt}
			\includegraphics[trim=3.9cm 9.5cm 4.8cm 9.7cm,clip=true,width=\columnwidth]{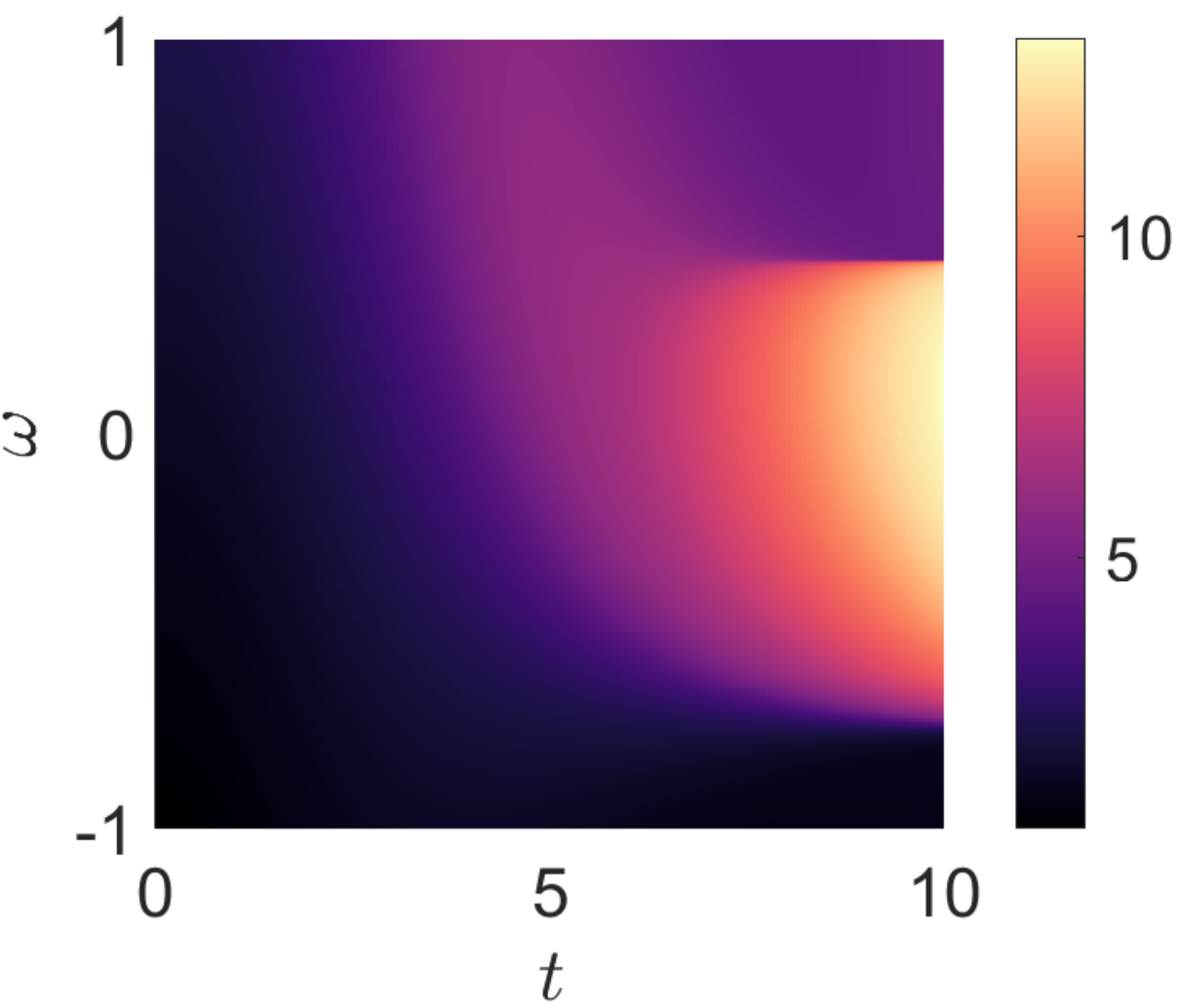}
		\end{subfigure}
		\begin{subfigure}[b]{0.43\columnwidth}
			\caption{$x_2(t,\omega)$}
			\vspace{-6pt}
			\includegraphics[trim=3.9cm 9.5cm 4.8cm 9.7cm,clip=true,width=\columnwidth]{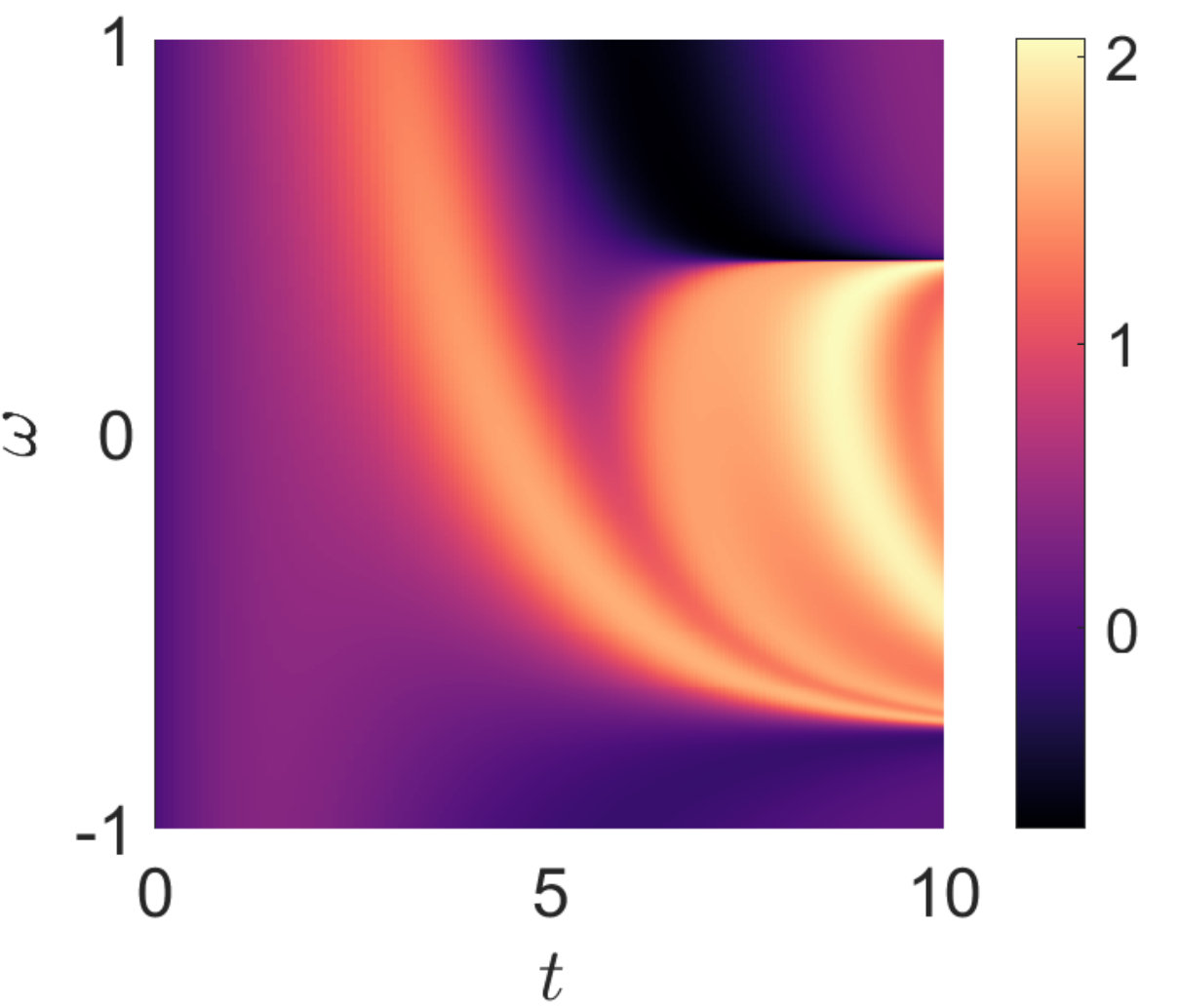}
		\end{subfigure}
		\begin{subfigure}[b]{0.43\columnwidth}
			\caption{$x_{(d),1}(t,\omega)$}
			\vspace{-6pt}
			\includegraphics[trim=3.9cm 9.5cm 4.8cm 9.7cm,clip=true,width=\columnwidth]{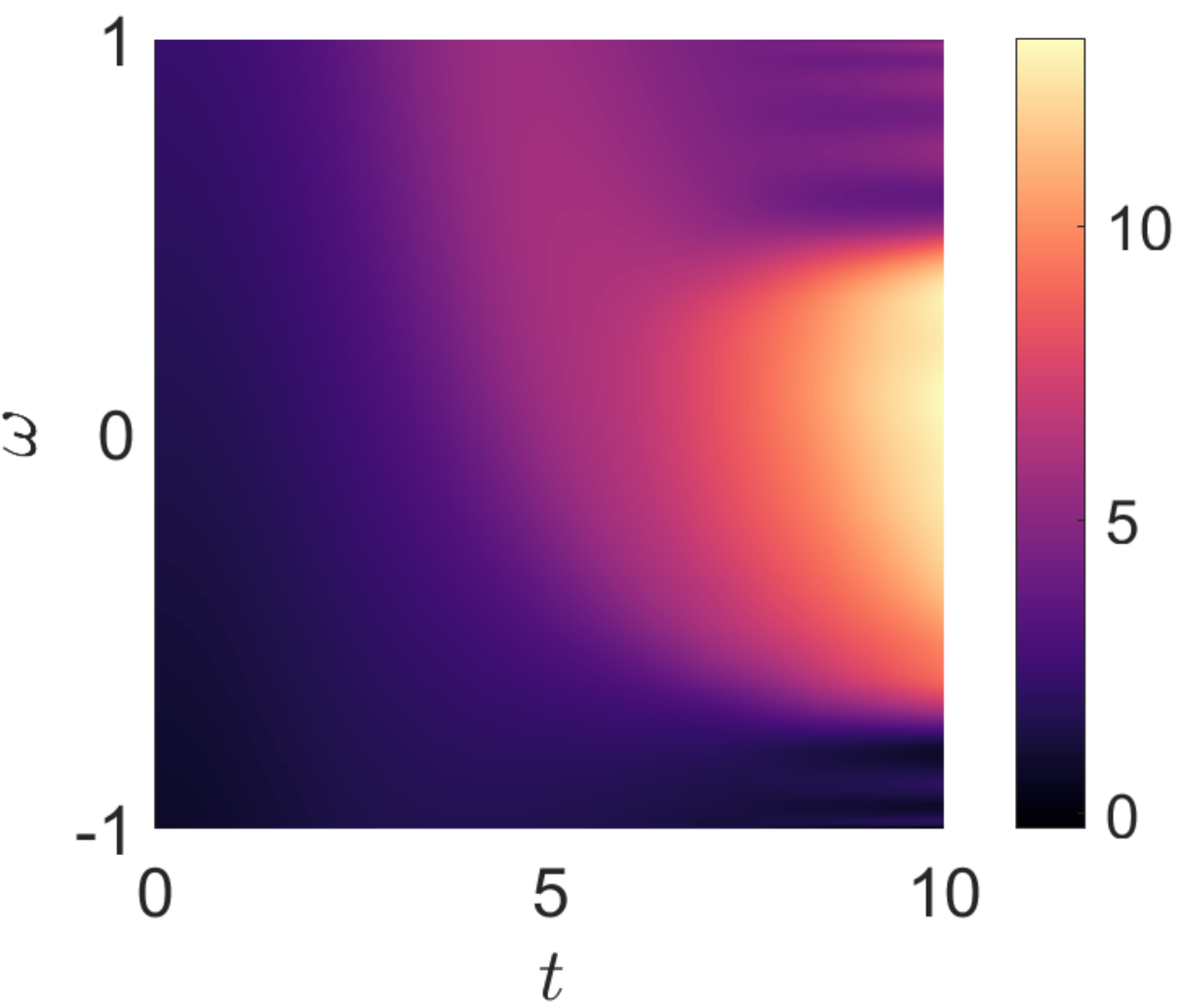}
		\end{subfigure}
		\begin{subfigure}[b]{0.43\columnwidth}
			\caption{$x_{(d),2}(t,\omega)$}
			\vspace{-6pt}
			\includegraphics[trim=3.9cm 9.5cm 4.8cm 9.7cm,clip=true,width=\columnwidth]{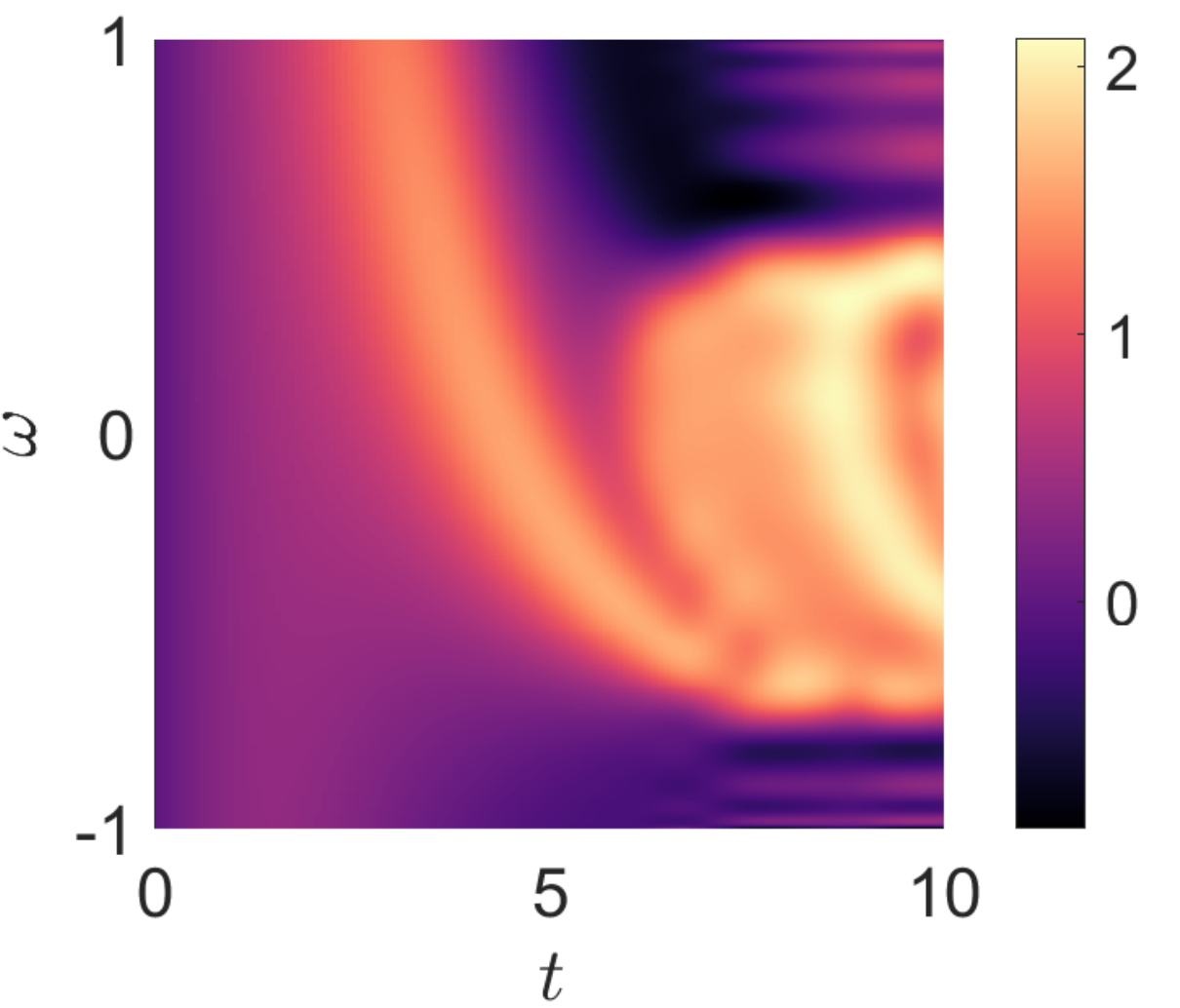}
		\end{subfigure}
		\caption{Comparison of Monte Carlo ($N=500$) and gPC surface approximations using the projection $(d=20,q=21)$  and coupled coefficient method for scenario $1$.}
		\label{fig:surfaces}
		\vspace*{-.5cm}
	\end{figure}

	It should be recognized, that the framework struggles to capture the bifurcation, seen in Fig. \ref{fig:surfaces}. Here we visualized the surfaces, $x_1(t,\omega)$ and $x_2(t,\omega)$, and the gPC approximations for scenario 1. The magnitude of the high-order coefficients remains significant, compared to scenario $2$, where they gradually fade out. This is due to the high frequency content of the corresponding basis polynomials which the framework tries to employ to model the discontinuity. As a result, one can observe that the discontinuity `spreads out' and causes a ripple effect in the otherwise smooth regions.

	\begin{table*}[t]
		\centering
		\caption{Comparison between RMSE values between $x_1(t,\omega)$ and ${x}_{(d),1}(t,\omega)$ for varying $q$ and $d\leq q$, using the LS and the projection method. The left values correspond with the decoupled method, the right values with the coupled method.}
		\label{tab:rmse}
		\begin{tabular}[\textwidth]{l|c|cccc|cccc}
			\toprule
			\multicolumn{2}{c}{} & \multicolumn{4}{c}{scenario 1} & \multicolumn{4}{c}{scenario 2} \\
			\midrule
			& $q\backslash d$ & 2 & 4 & 10 & 20 & 2 & 4 & 10 & 20 \\
			\midrule
			\multirow{4}{*}{LS} & 3 & 19.8 / 19.8 & & & & 3.09E-2 / 3.09E-2 & & &  \\
			& 5 & 13.4 / 28.6 & 9.87 / 9.87 & & & 2.99E-2 / 1.57E-1 & 1.05E-3 / 1.05E-3 & &  \\
			& 11 & 13.5 / 24.1 & 5.07 / 7.43 & 2.54 / 2.54 & & 3.29E-2 / 1.77E-1 & 0.91E-3 / 1.84E-2 & 4.66E-6 / 4.66E-6 &  \\
			& 21 & 13.8 / 24.9 & 5.11 / 7.49 & 2.13 / 2.17 & 1.49 / 1.49 & 3.42E-2 / 1.81E-1 & 0.95E-3 / 1.94E-2 & 3.35E-6 / 3.35E-6  & 9.64E-8 / 9.64E-8  \\
			\midrule
			\multirow{4}{*}{PM}	& 3 & 19.8 / 19.8 & & & & 3.09E-2 / 3.09E-2 & & & \\
			& 5 & 14.2 / 26.1 & 9.87 / 9.87 & & & 2.74E-2 / 1.26E-1 & 1.05E-3 / 1.05E-3& & \\
			& 11 & 12.1 / 20.8 & 5.10 / 5.96 & 2.54 / 2.54& & 2.74E-2 / 1.32E-1& 0.80E-3 / 1.22E-2 & 4.66E-6 / 4.66E-6& \\
			& 21 & 12.1 / 22.2 & 5.06 / 6.86 & 2.16 / 2.23 & 1.49 / 1.49 & 2.74E-2 / 1.32E-1 & 0.80E-3 / 1.22E-2& 3.05E-6 / 3.05E-6 &  9.64E-8 / 9.64E-8  \\
			\bottomrule
		\end{tabular}
	\end{table*}

\begin{figure}[t!]
		\centering
		\begin{subfigure}[b]{\columnwidth}
			\centering
			\vspace{-7pt}
			\includegraphics[trim=4.5cm 9.25cm 5.5cm 9.75cm,clip=true,width=.45\columnwidth]{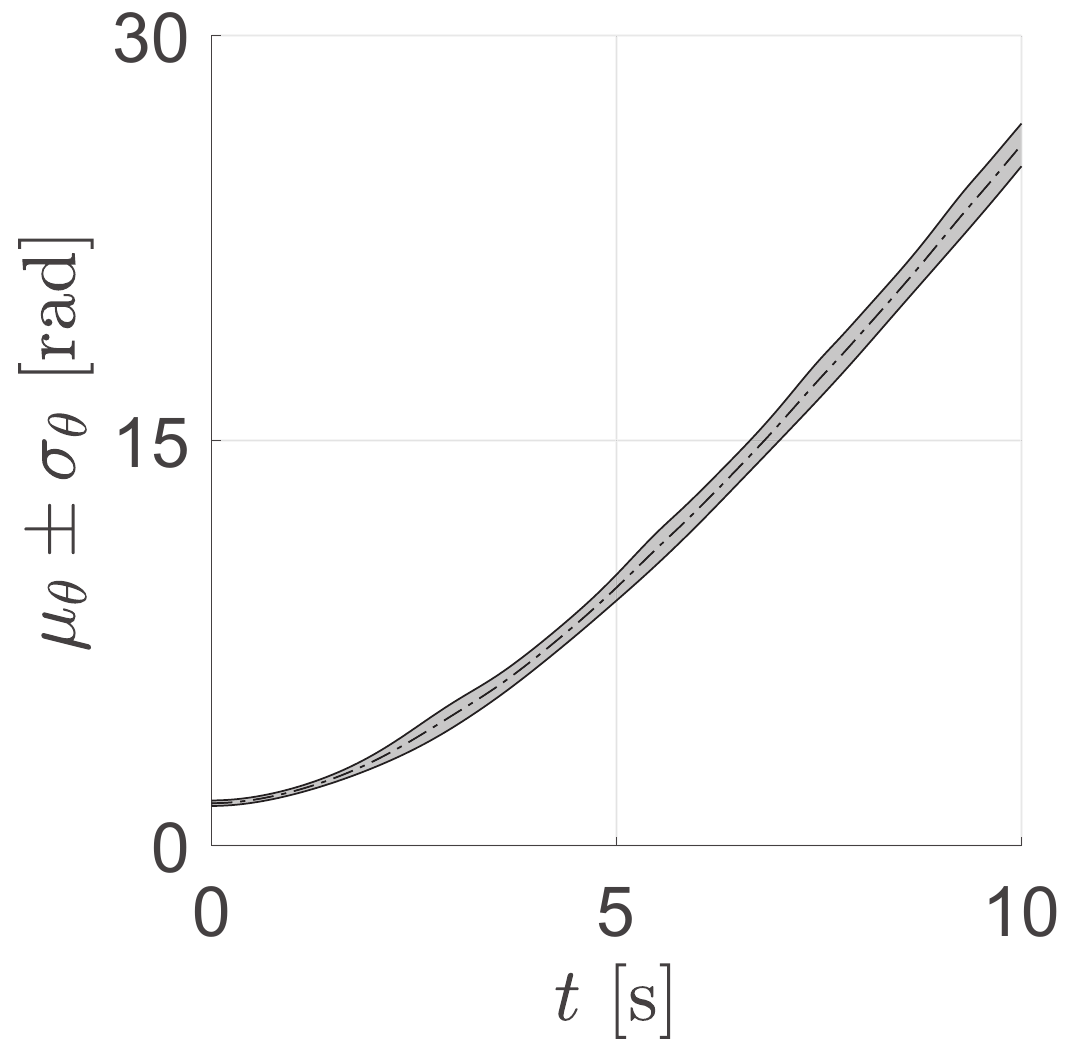}
			\includegraphics[trim=4.5cm 9.25cm 5.5cm 9.75cm,clip=true,width=.45\columnwidth]{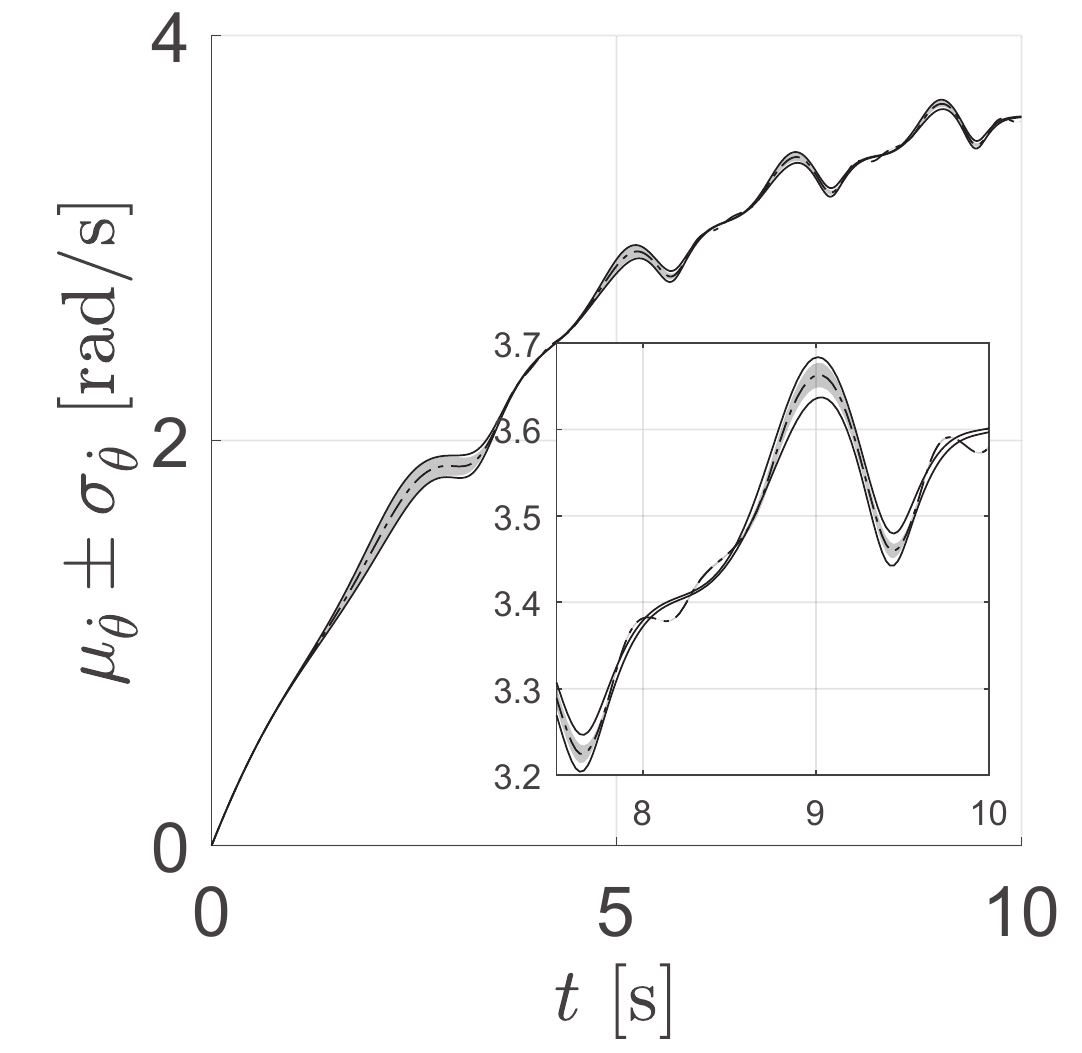}
			\caption{decoupled coefficient determination}
		\end{subfigure}
		\begin{subfigure}[b]{\columnwidth}
			\centering
			\includegraphics[trim=4.5cm 9.25cm 5.5cm 9.75cm,clip=true,width=.45\columnwidth]{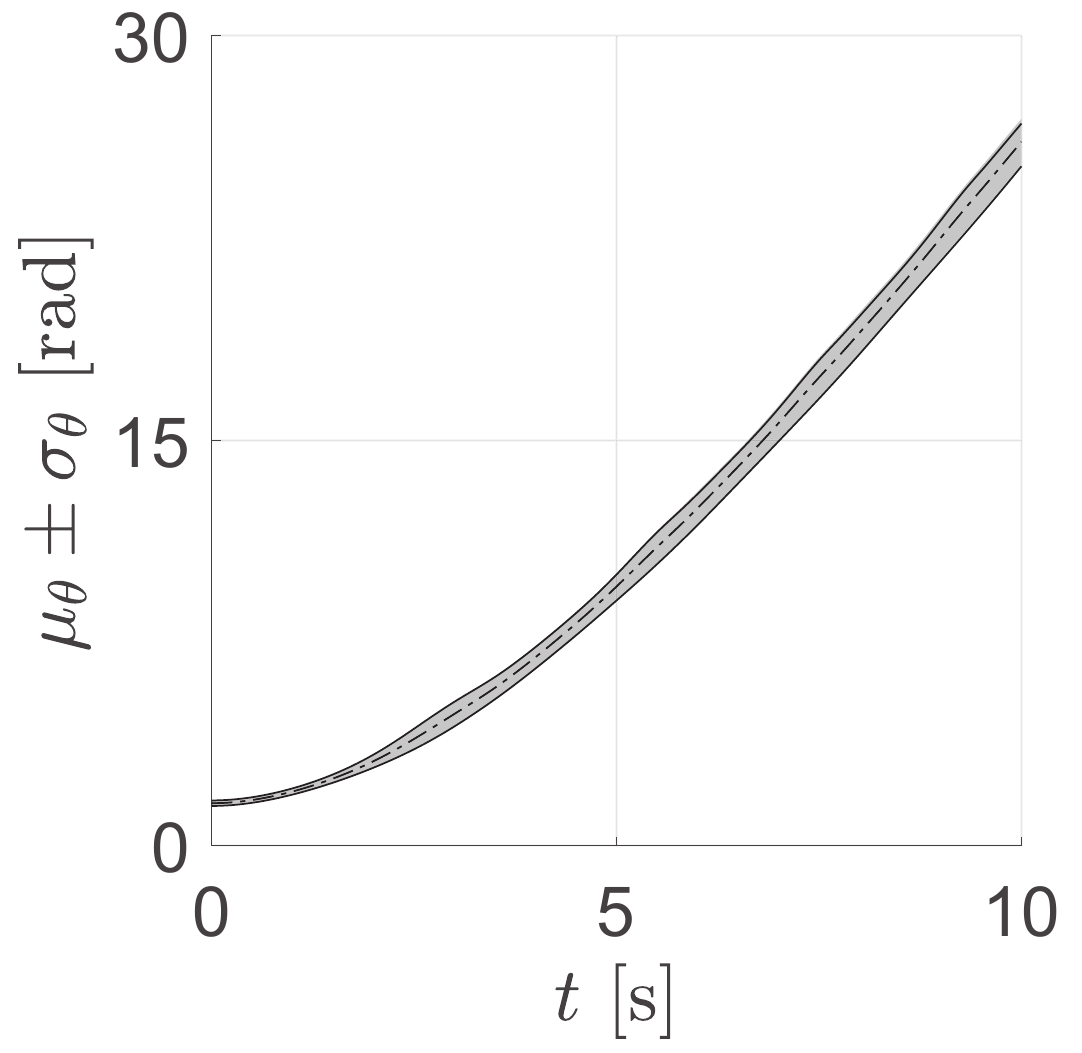}
			\includegraphics[trim=4.5cm 9.25cm 5.5cm 9.75cm,clip=true,width=.45\columnwidth]{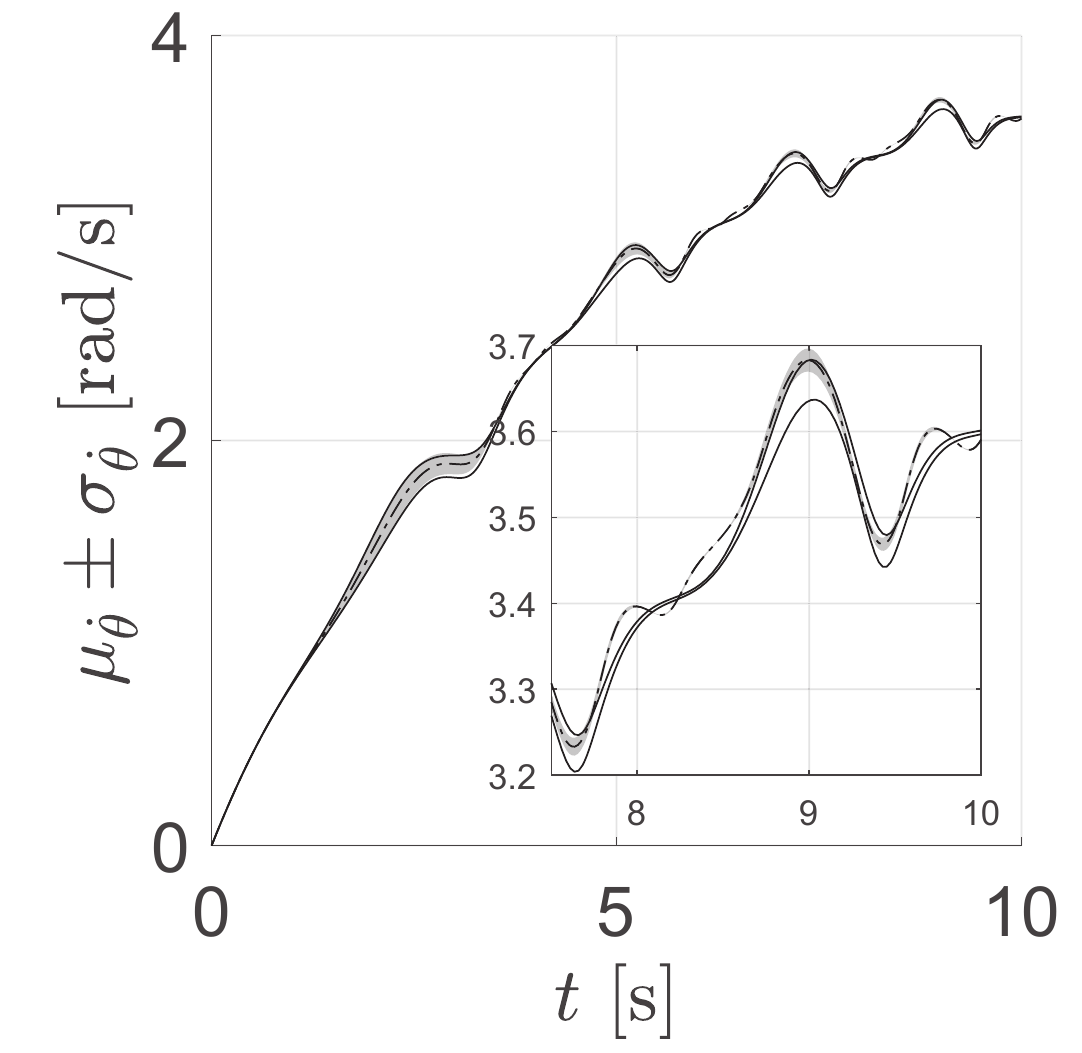}
			\caption{coupled coefficient determination}
		\end{subfigure}
		\caption{First two moments (grey shaded confidence interval) obtained with (\ref{eq:moments}) in scenario $2$ for both the decoupled as coupled coefficients $(d=2,q=5)$, compared with Monte Carlo $99\%$-confidence interval in black ($N=500$).}
		\label{fig:timeplots}
		\vspace*{-.5cm}
	\end{figure}

	When we stay outside the bifurcation range, as is the case in scenario $2$, the gPC framework proofs to be an efficient tool to quantify uncertainty. Such is illustrated by Fig. \ref{fig:timeplots}, where the mean value with $99\%$-confidence interval obtained with low chaos order, is compared to that obtained through $500$ MC simulations. Moreover, regarding the first two moments, experiments with more complicated input policies (not documented here) learned that, as long as  bifurcations are avoided, gPC behaves excellent and superior performance of either the decoupled or coupled dynamic coefficient determination approach is case dependent.
	
	\subsection{Robustified optimal start-up}
	In order to illustrate the importance of the robustification, we consider a smoothed ramp reference trajectory. The ramp is chosen to closely resemble the time evolution of $x_1(t)$ for $\theta(0) = \pi/2,\dot{\theta} = 0$ in scenario $1$, so to provoke a bifurcation. 
	\begin{equation}
	\label{eq:referencetraj}
	r_1(t) = \left(\tfrac{\pi}{8}t^2 + \tfrac{\pi}{2}\right)\cdot\mathbb{H}(2-t) + \tfrac{\pi}{2}t\cdot\mathbb{H}(t-2)
	\end{equation}
	
	We desire to obtain a control policy that avoids bifurcational behaviour regardless the actual value of $\theta_0$. The control horizon is fixed to $T = 10$ \si{\second} and the control discretisation to $N = 40$. The stochastic cost functional is defined as in (\ref{eq:costfunctional_stochastic}). The uncertainty is quantified using the Legendre polynomials, with a chaos order $d=7$. We used Gaussian quadrature with $q=15$. The PM method is combined with the decoupled dynamic coefficient determination. 
	
\begin{figure}[t!]
		\centering
		\begin{subfigure}[b]{0.43\columnwidth}
			\vspace{-7pt}
			\includegraphics[trim=4.5cm 9.25cm 5.5cm 9.75cm,clip=true,width=\columnwidth]{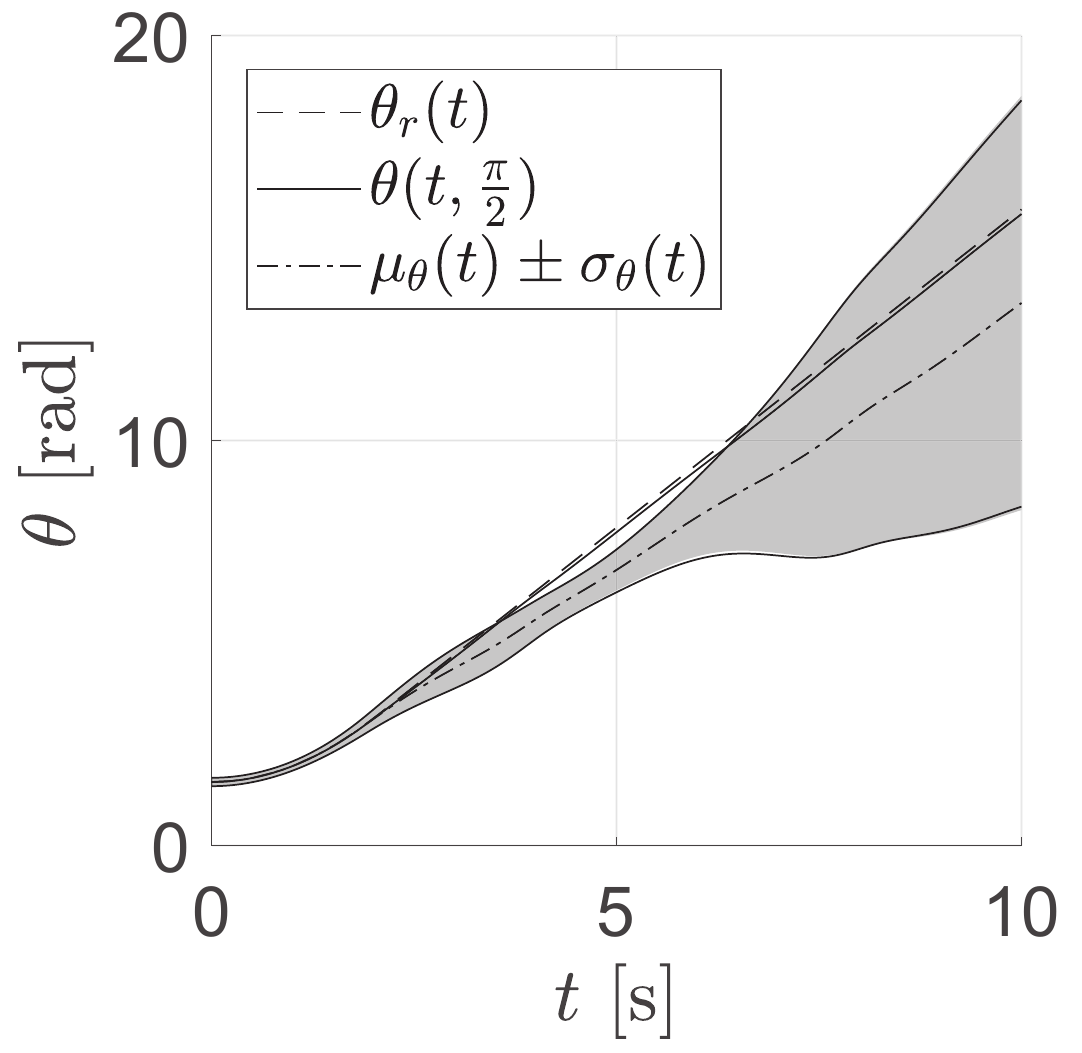}
		\end{subfigure}
		\begin{subfigure}[b]{0.43\columnwidth}
			\vspace{-7pt}
			\includegraphics[trim=4.5cm 9.25cm 5.5cm 9.75cm,clip=true,width=\columnwidth]{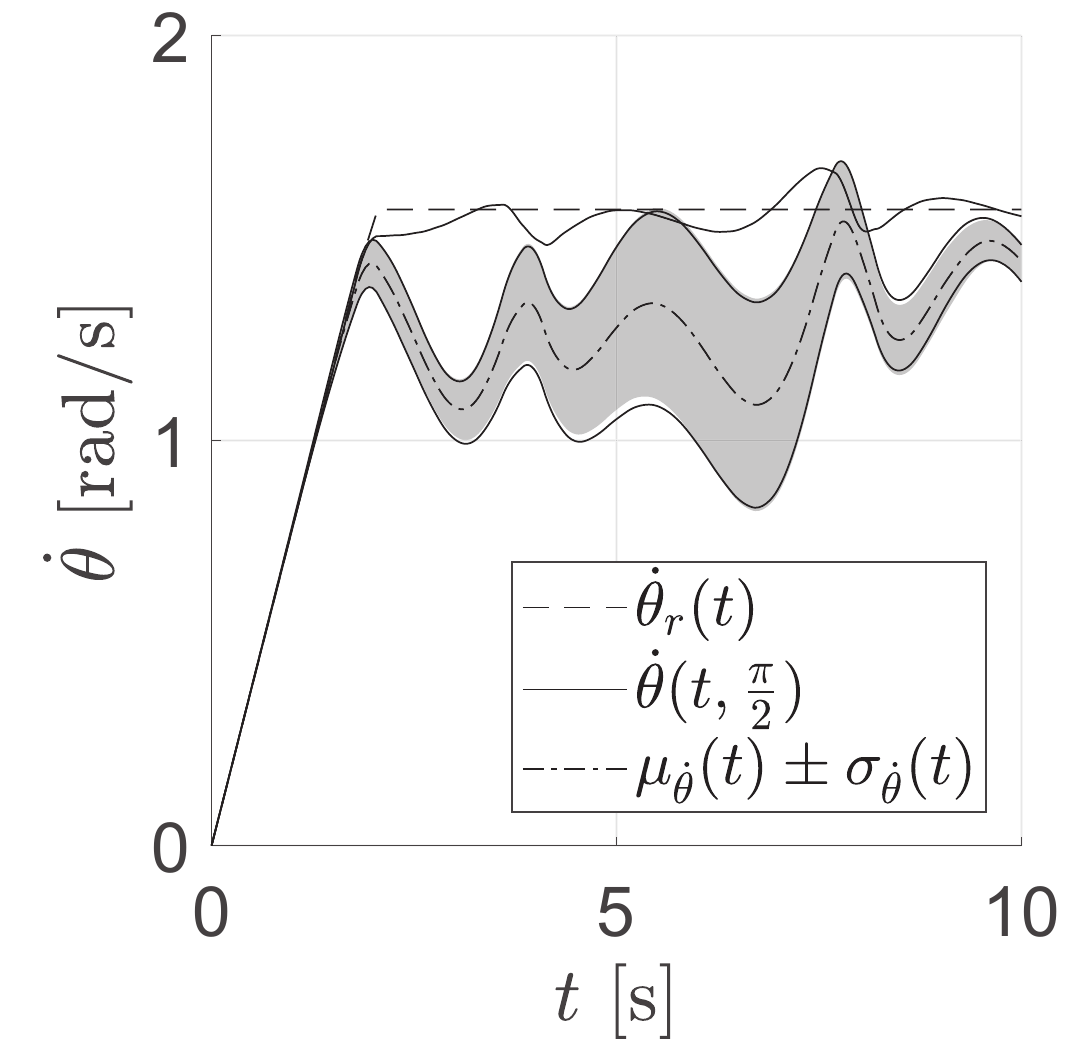}
		\end{subfigure}
		\caption{System behaviour for the computed torque control given reference trajectory (\ref{eq:referencetraj}). Visualization of the first two moments (grey shaded confidence interval) obtained with the PM decoupled method ($d=7,q=15$), compared to MC $99\%$-confidence interval (black). Also depicted are the reference trajectories and the actual trajectory for $\theta_0 = \frac{\pi}{2}$.}
		\label{fig:computedtorquecontrol}
	\end{figure}

\begin{figure}[t!]
		\centering
		\begin{subfigure}[b]{0.43\columnwidth}
			\vspace{-7pt}
			\includegraphics[trim=4.5cm 9.25cm 5.5cm 9.75cm,clip=true,width=\columnwidth]{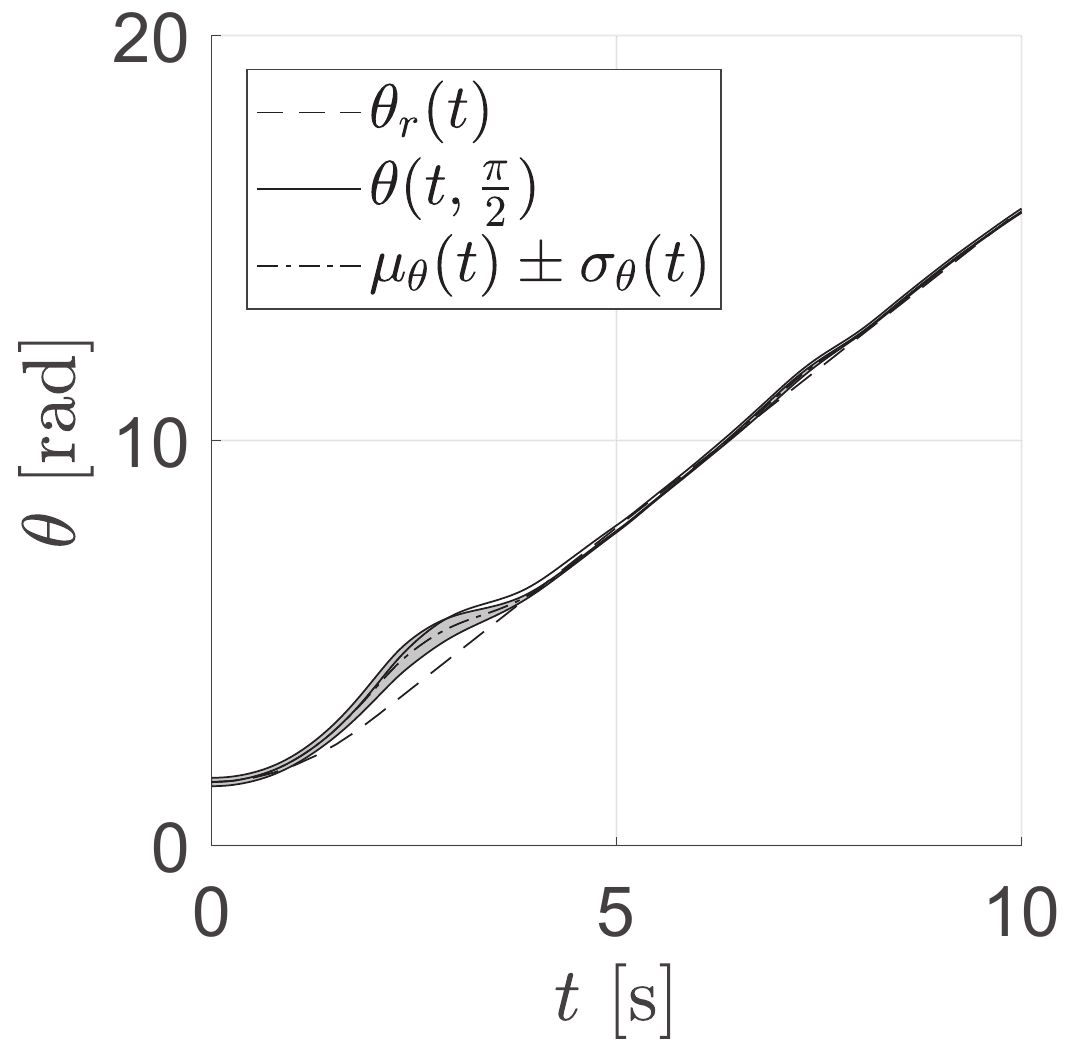}
		\end{subfigure}
		\begin{subfigure}[b]{0.43\columnwidth}
			\vspace{-7pt}
			\includegraphics[trim=4.5cm 9.25cm 5.5cm 9.75cm,clip=true,width=\columnwidth]{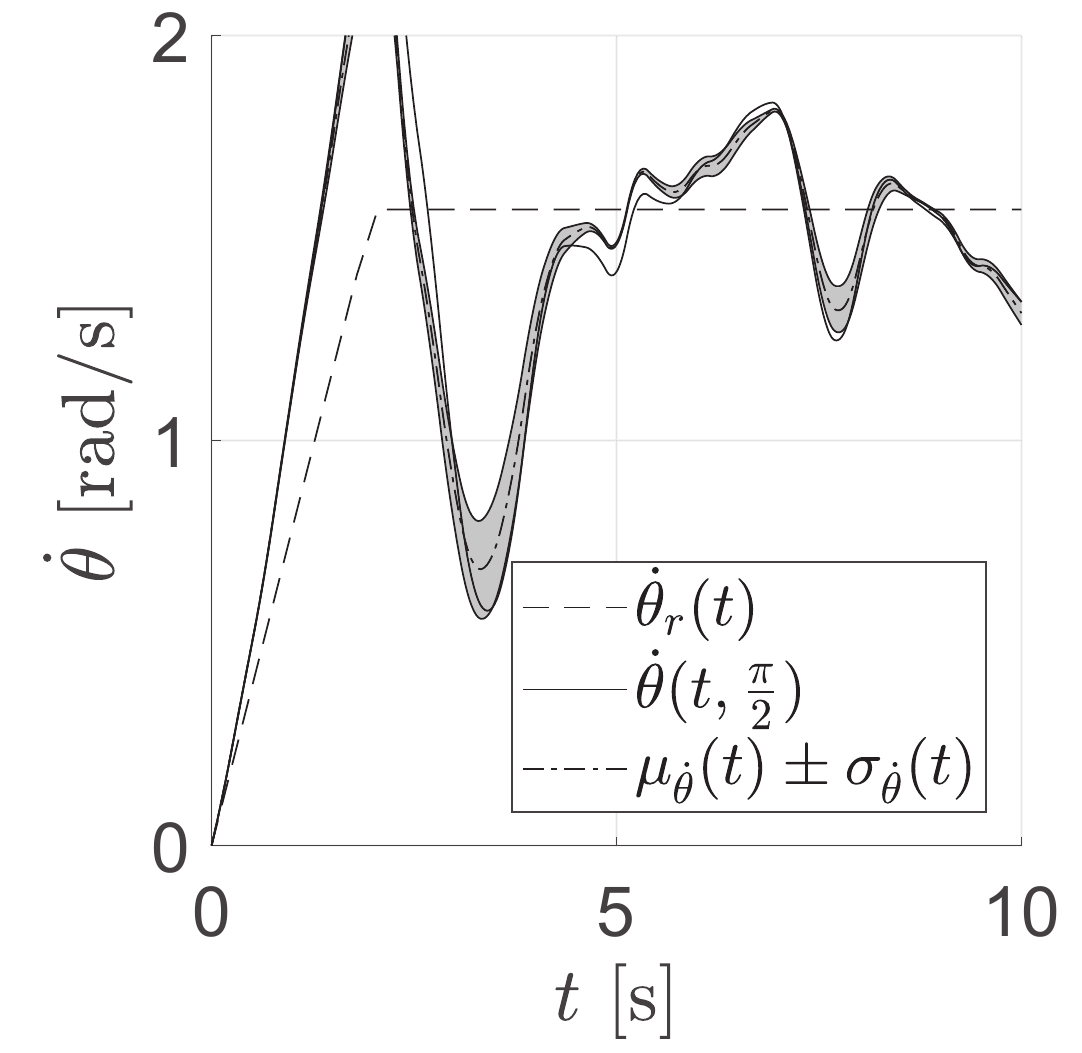}
		\end{subfigure}
		\caption{System output with a robustified control policy. The colours and the signals presented are defined as in Fig. \ref{fig:computedtorquecontrol}.}
		\label{fig:computedtorquecontrol_NMBC}
		\vspace*{-.5cm}
	\end{figure}

	Fig. \ref{fig:computedtorquecontrol} depicts the system behaviour as obtained for a forward computed torque control strategy. One can observe the rapidly increasing confidence interval implying that, as anticipated, bifurcational behaviour is present. The full black line depicts the system response in the nominal case, i.e. $\theta_0 = \frac{\pi}{2}$. The small discrepancy between the reference trajectory is a discretisation error resulting the coarse sampling of the actual torque signal. An illustrative stochastic cost functional was defined with ${\mathrm{Q}} = \left[\begin{smallmatrix}
	1 & 0 \\ 0 & 1
	\end{smallmatrix}\right]$, $\mathrm{R} = 1$ and $\epsilon=\tfrac{2}{5}$. System behaviour as obtained with the stochastic optimal control framework is presented in Fig. \ref{fig:computedtorquecontrol_NMBC}. The cost could be reduced from $K = 227$ for the computed torque control to $K = 21$ using the solution of ($\ref{eq:stochasticproblem}$). The stochastic optimal control framework successfully avoids the bifurcation behaviour, which implies that for every value of $\theta_0$ a tracking motion is achieved. Resulting the definition of ${\mathrm{Q}}$, we find that quite some oscillations are still present in $\mu_{\dot{\theta}}(t)$, yet remain acceptable with respect to $\mu_{\theta}(t)$.
	
	%
	
	\section{Conclusion}
	The gPC framework has proven to be an efficient mathematical tool to quantify uncertainty in general and to propagate parametric uncertainty in the context of dynamical systems. We have employed this framework as a mathematical tool to deal with uncertainty in the optimal control framework. To that end we introduced a stochastic reformulation of a quadratic tracking cost functional that penalizes performance and robustness independently. Introducing the general polynomial chaos framework into this specific problem formulation allowed to formulate an equivalent deterministic quadratic optimal control problem in function of the expansion coefficients. The method is tested numerically on the start-up behaviour of a realistic nonlinear drivetrain. The method proved to be tractable and successful and was able to avoid bifurcational behaviour compared to a feedforward computed torque control. 
	
	\section*{Acknowledgements}
	The autors are associated to EEDT Decision and Control, Flanders Make, Belgium and acknowledge support of FWO project G.0D93.16N and Flanders Make project EMODO.
	
	\bibliographystyle{unsrt}
	\bibliography{references}
	
	\appendices
	\section{Details on stochastic OC reformulation}
	\label{appendix:details_PCeSOC}
	Evaluation of the expected value of the original cost functional can be decomposed into two distinct contributions, i.e. $\mathbb{E}\left[J\right] = (\mathrm{I}) + (\mathrm{II})$. We elaborate each seperatly.
	
	The integrand of the first contribution is given by
	\begin{equation}
	\label{eq:expected}
	\begin{aligned}
	(\mathrm{I}) &= \int_{\Gamma} \int_{\mathcal{T}} (\smallsum_i \tilde{\mathbf{x}}_i\Psi_i - \mathbf{r})^\top{\mathrm{Q}}(\smallsum_j \tilde{\mathbf{x}}_j\Psi_j - \mathbf{r}) \text{d}\tau \rho\text{d}\uline{\omega}\\
	& = \int_{\mathcal{T}} \textstyle \sum_{i} \langle\Psi_i^2\rangle \|\tilde{\mathbf{x}}_i\|^2_\mathrm{Q} + \|\mathbf{r}\|^2_\mathrm{Q} - 2 \mathbf{r}^\top\mathrm{Q}\tilde{\mathbf{x}}_1 \text{d}\tau\\
	& = \int_\mathcal{T} \|{\tilde{\mathbf{X}}}\|^2_{\mathrm{D}\otimes\mathrm{Q}} + \|{\tilde{\mathbf{R}}}\|^2_{\mathrm{D}\otimes\mathrm{Q}} - 2{\tilde{\mathbf{R}}}^\top({\mathrm{D}}\otimes{\mathrm{Q}}){\tilde{\mathbf{X}}} \text{d}\tau\\
	&= \int_\mathcal{T}\|{\tilde{\mathbf{X}}}-{\tilde{\mathbf{R}}}\|^2_{\mathrm{D}\otimes\mathrm{Q}} \text{d}\tau
	\end{aligned}
	\end{equation}
	
	The second contribution can be evaluated as
	\begin{equation}
	\begin{aligned}
	(\mathrm{II}) &= \smallsum_i \smallsum_j \mathbf{u}_i^\top{\mathrm{R}}\mathbf{u}_j \int_\mathcal{T} D^{i}\mathrm{tri}(\tfrac{\tau}{\Delta}) D^{j}\mathrm{tri}(\tfrac{\tau}{\Delta})\text{d}\tau \\
	& =  \smallsum_i\smallsum_j \mathbf{u}_i^\top{\mathrm{R}}\mathbf{u}_j \Delta_{ij} = \|{\mathbf{U}}\|^2_{\mathrm{M}\otimes\mathrm{R}}
	\end{aligned}
	\end{equation}
	where $\Delta_{ij}$ correspond with the element of the matrix $\mathrm{M}$.
	
	The Frobenius norm of the covariance can be evaluated as
	\begin{equation}
	\begin{aligned}
	\label{eq:covariance}
	\|\mathrm{cov}[\mathbf{x}]\|_F^2 &\approx \mathrm{tr}\left(\smallsum_i\langle\Psi_i^2\rangle \tilde{\mathbf{x}}_i \tilde{\mathbf{x}}_i^\top - \tilde{\mathbf{x}}_1 \tilde{\mathbf{x}}_1^\top\right) \\
	&= \smallsum_i \smallsum_j \langle\Psi_i^2\rangle \tilde{x}_{ij}\tilde{x}_{ij} - \smallsum_j \tilde{x}_{1j}\tilde{x}_{1j} \\
	&=  \smallsum_i \langle\Psi_i^2\rangle \|\tilde{\mathbf{x}}_i\|^2_2 - \|\tilde{\mathbf{x}}_1\|^2_2 = \|{\tilde{\mathbf{X}}}\|^2_{\mathrm{E}\otimes\mathrm{I}}
	\end{aligned}
	\end{equation}
	
	In conclusion, note that one can combine the integrand in equation (\ref{eq:expected}) and equation (\ref{eq:covariance}) to yield an alternative interpretation of the stochastic problem formulation in (\ref{eq:stochasticproblem})
	\begin{equation}
	\begin{multlined}
	\epsilon(\ref{eq:expected}) + \left(1-\epsilon\right)(\ref{eq:stochasticproblem}) \\ 
	\begin{aligned}
	&= \epsilon\|{\tilde{\mathbf{X}}}-{\tilde{\mathbf{R}}}\|^2_{\mathrm{D}\otimes\mathrm{Q}} + \left(1-\epsilon\right) \|{\tilde{\mathbf{X}}}\|^2_{\mathrm{E}\otimes\mathbb{I}}
	\\
	&= \epsilon\smallsum_{i=1}^P \langle\Psi_i^2\rangle\|\tilde{\mathbf{x}}_i - \delta_{1i}\mathbf{r}\|^2_{\mathrm{Q}} + \left(1-\epsilon\right)\smallsum_{i=2}^P \langle\Psi_i^2\rangle \|\tilde{\mathbf{x}}_i\|^2_2\\
	&= \|\tilde{\mathbf{x}}_1 - \mathbf{r}\|^2_{\epsilon\mathrm{Q}} + \smallsum_{i=2}^P \langle\Psi_i^2\rangle \|\tilde{\mathbf{x}}_i\|^2_{\epsilon\mathrm{Q} + \left(1-\epsilon\right){\mathrm{I}}} \\
	&\approx\|\mathbb{E}[\mathbf{x}]- \mathbf{r}\|^2_{\epsilon\mathrm{Q}} + \|\mathrm{cov}[{\mathrm{S}}\cdot\mathbf{x}]\|^2_F
	\end{aligned}
	\end{multlined}
	\end{equation}
	in the special case where $\mathrm{{S}}^2 = \epsilon\mathrm{{Q}} + \left(1-\epsilon\right){\mathrm{I}}$.

\end{document}